\documentclass[12pt]{amsart}

\date{\today}
\usepackage [utf8]{inputenc}
\usepackage[math]{iwona}
\usepackage{tikz,fullpage,multicol}
\usetikzlibrary{arrows,%
                petri,%
                topaths}%
\usepackage{tkz-graph,tkz-berge}
\usepackage{amsmath,amsfonts, amssymb, version,  enumerate,color,hyperref,graphicx}
\usepackage{layout}
\definecolor{lightgrey}{rgb}{0.3,0.3,0.3}

\newtheorem{Theorem}{Theorem}
\theoremstyle{change}

\newcommand\df[1]{{\emph{#1}}\index{#1}}

\newcommand\N{\mathbb{N}}
\newcommand\Z{\mathbb{Z}}

\newcommand\R{\mathbb{R}}
\newcommand\C{\mathbb{C}}

\renewcommand\L{\mathcal{L}}

 \newcommand\underlabel[2]{\underset{#1}{\underbrace{#2}}}

\newcommand\norm[1]{\lVert #1 \rVert}

\newcommand\summ[3]{\displaystyle\sum_{#1 = #2}^{#3}}

\newcommand\BC[2]{\binom {#1}{#2}}

\newtheorem{Ex}[Theorem]{Example}
\newtheorem{Def}[Theorem]{Definition}
\newtheorem{Lemma}[Theorem]{Lemma}
\newtheorem{Proposition}[Theorem]{Proposition}
\newtheorem{Ax}[Theorem]{Ax}
\newtheorem{Corollary}[Theorem]{Corollary}

\newtheorem{kNote}[Theorem]{Note}
\newtheorem{Nota}[Theorem]{Notation}
\newtheorem{Exer}[Theorem]{Exercise}
\newtheorem{Rem}[Theorem]{Remark}
\newtheorem{Fac}[Theorem]{Fact}
\newtheorem{kNotes}[Theorem]{Notes}

\newenvironment{Example}{\begin{Ex} \rm}{\end{Ex}\goodbreak}
\newenvironment{Definition}{\begin{Def}\rm}{\end{Def}\goodbreak}

\newenvironment{Proof}{\noindent\textbf{Proof. }}{\null\hfill$\square$\medskip}
   \newenvironment{Remark}{\begin{Rem}\rm}{\end{Rem}}

\newcommand\ismc\cong
\renewcommand\P{\mathcal{P}}
\newcommand\poly[1]{\P^{#1}}

\DeclareSymbolFont{extrasymbols}{OMS}{cmsy}{m}{n}
\DeclareMathDelimiter{\lVert}
  {\mathopen}{extrasymbols}{"6B}{largesymbols}{"0D}
\DeclareMathDelimiter{\rVert}
  {\mathclose}{extrasymbols}{"6B}{largesymbols}{"0D}
\renewcommand\norm[1]{\lVert #1 \rVert}

 \let\epsilon\varepsilon
 
\makeatletter
  \def\footnotestar{\xdef\@thefnmark{}\@footnotetext}
\makeatother

\title{Revisiting products and powers of $(m,p)$ and $(m,\infty)$-isometries}
\author{Michael Mackey}
\email{mackey@maths.ucd.ie}
\address{University College Dublin}

\begin{document}
\maketitle

\begin{abstract}
  We review known results concerning powers and
  products of $(m,p)$-isometries with a view to providing elementary
  proofs based
  on properties of polynomials.  We consider also the situation when
  $p=\infty$ where we find elements of graph theory and combinatorics
  arise naturally.
\end{abstract}

\section{Introduction}
An operator $T$ on a normed space $X$, $T:X\to X$, is termed an
$(m,p)$-isometry if
\begin{equation}\summ k0m \BC mk(-1)^k \norm{T^k x}^p =0
  \label{eq:mp}
\end{equation}
for all
$x\in X$.  Here $m\in \N$ and, usually, $p\in(0,\infty)$.   This definition was given originally in \cite{MR2859754},
generalising previous notions of $m$-isometry on Hilbert or Banach
space (cf. \cite{MR1037599,AglerStankus,MR936999,MR2654806}).  This has been a well
studied topic in
recent years and many fundamental facts have been established.   For
example, Bermúdez et al. in
\cite{MR2911496} show that a power of an $(m,p)$-isometry is
also an $(m,p)$-isometry and, in \cite{MR2993366}, that
if the commuting operators $S$ and $T$ are $(m,p)$ and
$(n,p)$-isometries respectively
then $ST$ is an $(m+n-1,p)$-isometry, strengthening a result
covering the case $n=2$ proven in \cite{MR2654806}.  
We will also look at whether these facts about powers of
$(m,p)$-isometries carry over to the
$p=\infty$ situation, using the concept of a $(m,\infty)$-isometry
given in \cite{MR2852193}.   While the answer is negative,
nevertheless we will see an interesting transformation of the setting
from linear operators to combinatorics of binary sequences.

The main results presented in this note are known and due principally
to Bermúdez et al. \cite{MR2911496,MR2993366,MR3947303}, Bayart
\cite{MR2859754} and Hoffmann et al. \cite{MR2852193}.  Our motivation is
only to provide alternative and more general proofs of these facts
using elementary properties of polynomials which allow the results for operators
to be seen as an important application.

The approach is based on a fact established in
\cite{MR2852193} which gives an alternative characterisation of
$(m,p)$-isometry.    We denote the one-variable
polynomials over $\R$ of degree less than or equal to $m$ by $\poly m$.

\begin{Theorem}[{\cite[Remark 3.6]{MR2852193}}]\label{thm:poly}
 An operator $T:X\to X$ is an $(m, p)$-isometry if and only if, for
 every $x\in X$ there is a polynomial $q_x\in \poly{m-1}$ such
 that $\norm{T^n x}^p=q_x(n)$ for all $n\in\N^0$.
\end{Theorem}

It can be seen from either characterisation that the inverse of
an invertible $(m,p)$-isometry is also an $(m,p)$-isometry, but the
most direct route is to  apply \eqref{eq:mp} with $x$ replaced by
$T^{-m}x$.   The fact
that an $(m,p)$-isometry is an $(m+1,p)$-isometry is more readily
accessed via Theorem~\ref{thm:poly} since $\poly{m-1}\subset \poly m$.

\section{Powers of $(m,p)$-isometries}

We begin with a short alternative proof of the following.
\begin{Theorem}[{\cite[Thm 3.1]{MR2911496}}]\label{thm:pow}
  Let $X$ be a normed space, $T\in L(X)$ an $(m,p)$-isometry for
  $m\in\N, p\ge 1$.  Then $T^r$ is also an $(m,p)$-isometry for $r\in\N$.
\end{Theorem}

\begin{Proof}
  Theorem~\ref{thm:poly} provides a polynomial $q_x\in \P^{m-1}$ such
  that $q_x(n)=\norm{T^n x}^p$ for all $x$.  Define another polynomial $s_x\in \P^{m-1}$ by $s_x(t)=q_x(rt)$.  Then
  $\norm{(T^r)^nx}^p =\norm{T^{rn}x}^p =q_x(rn) = s_x(n)$ for all $n$ and
  so, again by Theorem~\ref{thm:poly}, $T^r$ is an $(m,p)$-isometry.
\end{Proof}

Suppose $T$ is an $(m,p)$-isometry, so that for $x\in X$ there exists $q_x\in
\poly{m-1}$ such that $q_x(n)=\norm{T^nx}^p$ for all $n\in\N^0$.  If
$T$ is invertible then
it follows $q_x(-d)=\norm{T^{-d}x}^p$ for $d\in\N$.  This is because
for all $n\ge d$ we have
$q_x(n-d)=q_{T^{-d}x}(n)$ and the fundamental theorem of algebra
extends this identity beyond these values of $n$.  In
particular, at $n=0$, we have
$q_x(-d)=q_{T^{-d}x}(0)=\norm{T^{-d}x}^p$.   Moreover, since the
values of $q_x$ are non-negative on $\Z$, $q_x$ cannot have odd
degree and this observation yields:

\begin{Proposition}[{\cite[Proposition 2.4(b)]{MR2859754}}]
  If $T\in L(X)$ is an invertible $(m,p)$-isometry and $m$ is even
  then $T$ is an $(m-1,p)$-isometry.
\end{Proposition}

Another fact that will be useful later and uses similar reasoning is
the following.
\begin{Lemma}\label{lem:min}
  Let $r,s \in \N$.  If $T^r$ is an $(m,p)$-isometry and $T^s$ is an $(n,p)$-isometry,
  then both are $(l,p)$-isometries where $l=\min(m,n)$.
\end{Lemma}

\begin{Proof}
  Let $q^r_x$ and $q^s_x$ be polynomials guaranteed by
  Theorem~\ref{thm:poly},  that is $q^r_x(n)=\norm{(T^r)^nx}^p$ and
  $q^s_x(n)=\norm{(T^s)^nx}^p$ have degrees bounded by $m-1$
  and $n-1$ respectively.   Since $q^r_x(ns)=q^s_x(nr)$ for all $n$,
  the fundamental theorem of algebra guarantees $q^r_x(sy)=q^s_x(ry)$
  for all $y\in\R$.  This implies the two polynomials have the same
  degree which cannot be greater than $\min(m-1,n-1)$.   Another
  appeal to
  Theorem~\ref{thm:poly} completes the proof.
\end{Proof}

Concerning a converse of Theorem~\ref{thm:pow}, let us remark that ``roots'' of
$(m,p)$-isometries do not generally belong to the class, that is
$T^2$ may be a $(m,p)$-isometry while $T$ is not an $(m,p)$-isometry,
nor a $(\mu,p)$-isometry for any $\mu$.    A simple example is
provided by $T(x,y)=(2y,x/2)$ on $\C^2$.  Evidently, $T^2$ is an
(invertible) isometry
and so a $(1,p)$-isometry for all $p$
but $T$ is not an isometry, nor can it be a $(m,p)$-isometry for any
$m$ since $(\norm{T^n(0,1)})_n$ is bounded and non-constant, meaning
it cannot be interpolated by a polynomial.

The main positive result in this direction of ``root extraction'' is
provided by Bermúdez at al. in \cite[Theorem 3.6]{MR2993366} which
says, for example, that if $T^2$ and $T^3$ are both $(m,p)$-isometries
then $T$ is an $(m,p)$-isometry.    We provide an alternative proof of
this result, which in light of Lemma~\ref{lem:min} above, we may
restate as follows:

\begin{Theorem}\label{thm:root}
  Suppose both $T^r$ and $T^s$ are $(m,p)$-isometries, $r,s \in \N$, and let  $d=\gcd(r,s)$.
  Then $T^d$ is an $(m,p)$-isometry.
\end{Theorem}

\begin{Proof}
  For invertible $T$, a very short proof is available after the next
  section (see Remark~\ref{rem:sp}) but we will argue without that
  assumption.  There is no loss in generality in assuming
  $d=\gcd(r,s)=1$, for otherwise we may apply the $d=1$ case to
  $S^{r'}$ and $S^{s'}$ where $S=T^{d}$ and $r'=r/d, s'=s/d$.

  Theorem~\ref{thm:poly} provides us with interpolating polynomials
  $q^r_x$ and $q^s_x\in\poly{m-1}$ for $T^r$ and $T^s$ at $x$ respectively.
  Evidently, $q^r(ns)=q^s(nr)$ for all $n\in \N^0$ which, by the
  fundamental theorem of algebra, guarantees $q^r(sy)=q^s(ry)$ for all
  $y\in\R$.

  Define another polynomial $q\in\poly{m-1}$ by
  $q(y):=q^r(y/r)=q^s(y/s)$.  We will show that $q$ interpolates $T$
  at $x$.

  Since $\gcd(r,s)=1$ we may choose $a,b\in\Z$ with $ar+bs=1$  Without loss
  of generality, we may assume $a>0$ and $b<0$.  Indeed this may be
  guaranteed since $(a+sn)r+(b-rn)s=1$ for all $n\in\N$.

  For any $k\in\N^0$ and $n\in\N$, we have $k(a+sn)r =k+ k(rn-b)s$ where all terms
  are non-negative and so are valid powers of $T$.   The polynomials
  $q^r_x$ and $q^s_{T^kx}$ can now be related as follows.
  \begin{align*}
    q^r_x(k(a+sn)) &= \norm{(T^r)^{k(a+sn)}x}^p \\
                   &= \norm{(T^s)^{k(rn-b)}T^kx}^p\\
    &= q^s_{T^kx}(k(rn-b))
  \end{align*}
  for all $n\in\N$.   The agreement of the two polynomials on this
  infinite sequence implies that \[q^r_x(k(a+sy)) = q^s_{T^kx}(k(ry-b))\] for all
  $y\in\R$.  In particular, taking $y=b/r$, we have \[
    q(k)=q^r_x(k/r)=q^r_x(k(a+s\frac br)) = q^s_{T^kx}(0) =\norm{T^kx}^p.\]
  Again by Theorem~\ref{thm:poly}, this means that $T$ is an $(m,p)$-isometry.
\end{Proof}

\begin{Remark}\label{rem:5}
  None of the above results depend essentially on the linearity or
  boundedness of the linear operator $T$.  One could define an
  \df{$m$-sequence} to be a any real sequence $(a_n)_{n\in\N^0}$ which
  satisfies $\summ k0m (-1)^k\BC mk  a_k=0$ or, equivalently, is interpolated by
  a polynomial $q\in\poly{m-1}$ in the sense that $a_k=q(k)$ for all
  $k\in\N^0$.  An $(m,p)$-isometry is then just an operator for which
  $(a_k)_k= (\norm{T^kx}^p)_k$ is an $m$-sequence.  The results above
  can be stated and proved in the  context of $m$-sequences.  For
  example, Theorem~\ref{thm:pow} is an exposition of the fact that the
  regular subsequence $(a_{rk})_k$ of an $m$-sequence is also an $m$-sequence.
  Nevertheless, it is interesting and instructive to see constructive
  examples of $(m,p)$-operators as this is where
  applications may arise.  While it is easy to write down
  $m$-sequences, it is much more difficult to construct a linear
  operator which realises any such sequence.
\end{Remark}

\section{Polynomial matrices}

  The two-variable polynomials of degree $\le n$ in $x$ and
degree $\le m$ in $y$ will be denoted $\P^{m,n}[y,x]$ and this is a
subspace of the set of all two variable polynomials of degree
$\le n+m$.  (Our less standard ordering of the indeterminants here is
for later consistency with matrix notation.) For $r\in \P^{m,n}[y,x]$,
we have $r^x \in \P^m[y]$ and $r_y\in\P^n[x]$ on fixing a coordinate-
that is, $r^x(y)=r(y,x)=r_y(x)$.


\begin{Definition}
  An infinite matrix $A=(a_{ij})_{i,j\in\N}$ will be called an
    $(n,m)$-polynomial matrix if each
  row is interpolated by a polynomial of degree at most $n$ on $\N$,
  and each column is interpolated by a polynomial of degree at most
  $m$ on $\N$.
  That is, for each $i\in \N$, there exists $p_i\in \P^n[x]$
  such that $a_{ij}= p_i(j)$ and for each $j\in \N$, there exists
  $q_j\in \P^m[y]$ such that $a_{ij}=q_j(i)$ for all $i$.
\end{Definition}

\begin{Example}\label{ex:}
  The matrix \[A=
    \begin{bmatrix}
      5 &8 &11 &14 &17 &\cdots \\
      6 &7 &6 &3 &-2 &\cdots \\
      7 &6 &1 &-8 &-21 &\cdots \\
      8 &5 &-4 &-19 &-40 &\cdots \\
      \vdots & \vdots &\vdots &\vdots &\vdots &\ddots
    \end{bmatrix}\]
  is a $(2,1)$-polynomial matrix.
\end{Example}

It is easy to see that if $r\in\P^{m,n}[y,x]$ then
$a_{ij}:=r(i,j)$ is an $(n,m)$-polynomial matrix.   (Example~\ref{ex:}  above
is the matrix so generated by $r(y,x)=1+2x+x^2+y+xy-x^2y$.)
The following proposition provides the converse to this statement.

\begin{Proposition}\label{prop:matrix}
  Let $A=(a_{ij})_{i,j\in \N}$ be an $(n,m)$-polynomial matrix.
  Then there exists a two-variable polynomial $r\in\P^{m,n}[y,x]$ such
  that $r(i,j)=a_{ij}$ for all $i,j \in\N$.
\end{Proposition}

\begin{Proof}
  We have polynomials $p_1, p_2,\ldots \in\P^n[x]$ such that
  $p_i(j)=a_{ij}$ for all $j\in\N$ and polynomials $q_1,q_2,\ldots
  \in\P^m[y]$ such that $q_j(i)=a_{ij}$ for all $i$.

 Consider the data set $(1,p_1),
 (2,p_2),\ldots,(m+1,p_{m+1})$.  Briefly treating $p_k$ as a scalar, we interpolate this data set with the
Lagrange polynomial in the usual way, $r(y):= \summ k1{m+1} p_k
\frac{L_k(y)}{L_k(y_k)}$ where 
$L(y)=(y-1)(y-2)\cdots(y-({m+1}))\in\P^{m+1}[y]$ and $L_k(y)=\frac{L(y)}{y-k}\in \P^m[y]$.

Now recalling each $p_k$ as an element of $\P^n[x]$,  we gain a two-variable polynomial $r\in\P^{m,n}[y,x]$ via
\[r(y,x)=\summ k1{m+1} p_k(x)\frac{L_k(y)}{L_k(k)}.\]
By construction, this $r\in\P^{m,n}[y,x]$ interpolates the
first $m+1$ rows of $A$.  Indeed, noting that $L_k(i)=0$ for
$k\in\{1,2,\ldots,i-1,i+1,\ldots m+1\}$, we have
\begin{equation}
  r(i,j)=\summ k1{m+1}
  p_k(j)\frac{L_k(i)}{L_k(k)}=p_i(j)=a_{ij} \label{eq:star}
\end{equation}
for $i=1,\ldots,m+1$ and $j\in \N$.

We must show that $r$ interpolates all of $A$, not just these rows.
  (One might recall at this point the general failure of the fundamental
theorem of algebra for two variable polynomials.)  

For  $i>m+1$, we know $a_{ij}= q_j(i)$.
However \eqref{eq:star} implies that for all $k\in\{1,\ldots,m+1\}$, we
have $q_j(k)=a_{kj}=r(k,j)=r^j(k)$ and so, via the fundamental
theorem of algebra, as both $q_j$ and $r^j$
are of degree at most $m$ and agree at these $m+1$ distinct points,
they must agree at all points.   Thus $r(i,j)= r^j(i)=q_j(i)=a_{ij}$ for \emph{all} $i,j\in \N$
as required.
\end{Proof}

\begin{Corollary}\label{cor:diag}
  If $A$ is an $(n,m)$-polynomial matrix then the diagonal entries are
  interpolated by a polynomial of degree at most $m+n$ on $\N$.  That
  is, there exists $q\in\P^{m+n}$ such that $a_{kk}=q(k)$ for all $k$.
\end{Corollary}

\begin{Proof}
  Proposition~\ref{prop:matrix} guarantees that $a_{kl}=r(k,l)$ for
  some $r\in\P^{n,m}[x,y]$.   Defining $q\in\P^{n+m}$ by
  $q(\xi)=r(\xi,\xi)$ gives the statement.
\end{Proof}

\section{Products of $(m,p)$-isometries}

\begin{Theorem}[{\cite[Thm 3.3]{MR2993366}}]\label{thm:prod}  Let $T$ and $S$ be commuting operators which are $(m,p)$ and
  $(n,p)$-isometries respectively.  Then $TS$ is an
  $(m+n-1,p)$-isometry.
\end{Theorem}

Our proof is again quite short given the facts established above for
polynomial matrices.

\begin{Proof}
  The hypothesis provides that the matrix $A$ given by
  $a_{ij}=\norm{S^iT^jx}^p$ is an $(n-1,m-1)$-polynomial matrix.  Indeed, Row $i$ is given by the sequence $(a_{ij})_j = (\norm{T^j
    (S^ix)}^p)_j $ which, by Theorem~\ref{thm:poly}, equals
  $(q_{S^ix}(j))_j$ for some $q_{S^ix}\in \P^{n-1}$.  Similarly,
  there is $\P^{m-1}$-polynomial interpolation of every column.

  Therefore,
  by Corollary~\ref{cor:diag}, there exists $q\in \P^{n+m-2}$ such that
  $q(k)=a_{kk}$.  As $a_{kk}=\norm{(TS)^kx}^p$, Theorem~\ref{thm:poly}
  again implies that $TS$
  is an $(n+m-1,p)$-isometry. 
\end{Proof}

\begin{Remark}\label{rem:sp}
  A simple proof of Theorem~\ref{thm:root}, in the case of  invertible
  $T$, can now be given.   Let $ar+bs=d$ for $a,b\in\Z$.  Since $T^r$
  and $T^s$ are $(m,p)$-isometries, 
  so too are their inverses $T^{-r}$ and $T^{-s}$ and, by
  Theorem~\ref{thm:pow}, any natural power of these.  From
  Theorem~\ref{thm:prod}, the product $(T^r)^a(T^s)^b=T^d$ is a
  $(2m-1,p)$-isometry.  However, since $T^r$ is an $(m,p)$-isometry, Lemma~\ref{lem:min} now implies that
  $T^d$ is an $(m,p)$-isomety.
\end{Remark}

\section{What about $(m,\infty)$-isometries?}

The following definition of an $(m,\infty)$-isometry was first given
in \cite{MR2852193}.  It comes from separating the positive and
negative terms in \eqref{eq:mp}, taking the $p$th root and letting $p$
tend to infinity.

\begin{Definition}
  \label{def}
An operator $T \in L(X)$ is called an \emph{$(m, \infty)$-isometry}, if
\begin{align}
	\max_{\substack{k=0,...,m \\ k \ \textrm{even}}}\norm{T^{k}x} 
	= \max_{\substack{k=0,...,m \\ k \ \textrm{odd}}}\norm{T^{k}x}, \
    \ \ \forall x \in X. \label{minf}
\end{align}
\end{Definition}
Just as for finite $p$, the inverse of an (invertible) $(m,\infty)$-isometry is
also an $(m,\infty)$-isometry and any $(m,\infty)$-isometry is an $(m+1,\infty)$
isometry.  Again linearity, and even boundedness, of $T$ are not 
essential features- \cite[Example~3.15]{PHthesis} provides an unbounded
operator satisfying \eqref{minf} for $m=2$.

Remark that we may replace $x$ by $T^rx$ (for $r\in \N^0$) in \eqref{minf} to
gain the fact that $T$ is an $(m,\infty)$-isometry if, and only if, \[\max_{\substack{k=r,...,m+r \\ k \ \textrm{even}}}\norm{T^{k}x} 
	= \max_{\substack{k=r,...,m+r \\ k \ \textrm{odd}}}\norm{T^{k}x}, \
    \ \ \forall x \in X, r\in \N^0.\]

As ventured in Remark~\ref{rem:5}, little of the literature on $(m,p)$
isometries relies greatly on
the theory of bounded linear operators and generally the
original context of linear operators can be substituted in favour of
the setting of sequences interpolated by polynomials.    This is also
true for $(m,\infty)$-isometries and we choose to take this approach
in the following.

\begin{Definition}
  \label{defseq}
A real sequence $a=(a_n)_{n\in \N^0} $ is called an \emph{$(m, \infty)$-sequence}, if
\begin{align}
	\max_{\substack{k=r,...,m+r \\ k \ \textrm{even}}} a_k 
	= \max_{\substack{k=r,...,m+r \\ k \ \textrm{odd}}}a_k \
    \label{minfseq}
\end{align}
for all $r\in\N^0$.
\end{Definition}

Clearly, such a sequence has a maximum and it attains this maximum at
least twice in any $m+1$ consecutive terms.  An operator $T$ is an
$(m,\infty)$-operator if and only if, for every $x$, the sequence
$a_x=(\norm{T^kx})_k$ is an $(m,\infty)$-sequence.  What we can prove
for $(m,\infty)$-sequences will have implications for the associated
operators.

A \emph{regular} subsequence of $(a_n)_n$ is a subsequence of the form
$(a_{rn})_n$ for some $r\in\N$. 
For an operator $T\in\L(X)$ and $x\in X$, we use $a^T_x$ to denote the sequence
$(a^T_x(n))_n$ given by
$a^T_x(n)=\norm{T^{n}(x)}$.   Notice that if $r\in \N$ then
$a^{T^r}_x(n)= a^T_x(rn)$.  In particular, the sequence generated by
$T^r$ is a regular subsequence of the sequence generated by $T$.

Forgoing treatment of constant sequences and isometries, we start with
the simple case of $(2,\infty)$-sequences (and operators).  These have a very simple structure.
Indeed, $(a_n)_{n\in\N^0}$ is a $(2,\infty)$-sequence iff
$a_1=\max(a_0,a_2)$, $a_2=\max(a_1, a_3)$, etc.    The first of these
equations shows $a_0\le a_1$ and $a_2\le a_1$ while the second implies
$a_1\le a_2$.    Thus $a_0\le a_1 = a_2 = a_3 \cdots$ and we have a
$(2,\infty)$-sequence if and only if it takes the form $(a_0, a_1,
a_1, a_1,\ldots)$ where $a_0\le a_1$.      Immediately then, we have the
following.
\begin{Lemma}\label{lem:13}
  A $(2,\infty)$-isometry is an isometry on its range.
\end{Lemma}
\begin{Lemma}
  Any regular subsequence of a $(2,\infty)$-sequence is also a
  $(2,\infty)$ sequence.
\end{Lemma}
\begin{Corollary}[{\cite[Theorem 3.2]{MR3947303}}]\label{cor:14}
  A power $T^r$ of a $(2,\infty)$-isometry $T$ is also a $(2,\infty)$-isometry.
\end{Corollary}

\subsection{$(3,\infty)$ binary sequences}

As mentioned above, any $(m,\infty)$ sequence has a maximum and attains it
repeatedly.  We can simplify consideration to binary sequences where
$1$ represents the maximum and $0$ represents any value which is not a
maximum.  It is the binary $(m,\infty)$ sequence that  we concentrate on, that is,
the possible patterns of maximum-attaining terms of the sequence.  For
example, we saw above that a $(2,\infty)$ sequence is either of the
form $(0,1,1,1,\ldots)$ or $(1,1,1,\ldots)$.

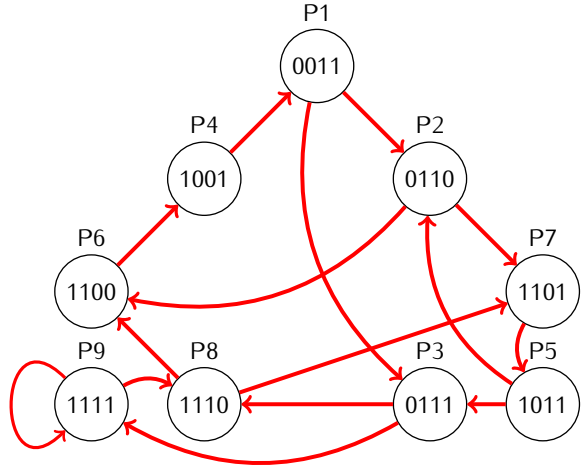
\begin{figure}[h]
 \centering
\begin{tikzpicture}[scale=0.75,transform shape]
   \SetUpEdge[lw         = 1.5pt,
            color      = red,
            labelcolor = white, labelstyle={above}]
            \SetVertexMath
            SetUpEdge[lw = 1.5pt, color = blue]
\GraphInit[vstyle=Normal] 
\SetGraphUnit{2}
\tikzset{VertexStyle/.append  style={fill}}
   \tikzset{EdgeStyle/.style={->}}

   \Vertex[style={label={P6}}]{1100}
   \NOEA[style={label={P4}}](1100){1001} \NOEA[style={label={P1}}](1001){0011} \SOEA[style={label={P2}}](0011){0110}\SOEA[style={label={P7}}](0110){1101}
   \SO[Lpos=-60,style={label={P5}}](1101){1011}\WE[style={label={P3}}](1011){0111} \SO[style={label={P9}}](1100){1111}\EA[style={label={P8}}](1111){1110}
   \Edges(1100,1001,0011,0110,1101)
   \Edges[style={bend right}](1101,1011)
   \Edges(1011,0111,1110,1101)
   \Edges[style={bend right}](0011,0111)
   \Edges[style={bend left}](0111,1111,1110)
   \Edges[style={bend left}](1011,0110,1100)
   \Edges(1110,1100)
   \Loop[dir=WE,dist=1.8cm,style={red,very thick}](1111)
    \end{tikzpicture}
 \caption{Directed graph for $(3,\infty)$ sequences}\label{fig:d1}
\end{figure}

The first, and indeed any consecutive, four terms of a $(3,\infty)$
sequence must have one of the following patterns:
\begin{center}
  \begin{multicols}{3}
    \begin{enumerate}[\bf P1.]
    \item 0,0,1,1
    \item 0,1,1,0
    \item 0,1,1,1
    \item 1,0,0,1
    \item 1,0,1,1
    \item 1,1,0,0
    \item 1,1,0,1
    \item 1,1,1,0
    \item 1,1,1,1
    \end{enumerate}
  \end{multicols}
\end{center}
Each of these may be considered as a vertex in a directed graph with
connections determined by possible patterns as we slide a window of
length 4 along the
$(3,\infty)$ sequence.  For example, {\bf P4} is connected to {\bf P1}
(we write {\bf P4}$\to${\bf P1}) and only {\bf P1}, while {\bf
  P1$\to$P2} and {\bf P1$\to$P3}.  The resulting graph is shown in Figure~\ref{fig:d1}.

Thus we represent $(3,\infty)$ sequences as infinite paths through
the graph.  Some of these paths are periodic, some are not.  Apart
from the cycle P9 $\to$  P9 representing a
constant sequence, the shortest cyclic paths are formed by the cycles
P1$\to$ P2 $\to$ P6 $\to$ P4 $\to$ P1 and  P7$\to$ P5 $\to$ P3$\to$ P8$\to$ P7.
For example, the former gives us the $(3,\infty)$ sequence
$(1,1,0,0,1,1,0,0,1,1,0,0,\ldots)$, or $(1100)$ for
short.  Notice however that if we take the regular subsequence given
by every second term then we have $(1,0,1,0,\ldots)$ or
$(10)$ and this is \emph{not} a $(3,\infty)$ sequence.

\begin{Corollary}\label{cor:fail}
  A power of an $(m,\infty)$-isometry is not necessarily an
  $(m,\infty)$-isometry.
\end{Corollary}

\begin{Proof}
  Consider the matrix $T$ giving the cyclic permutation of 
  coordinates on $\R^4$, $T:x=(x_1,x_2,x_3,x_4)\mapsto(x_2,x_3,x_4,x_1)$
  where the vector space is normed by \[\norm x=
  \max(2|x_1|,2|x_2|,|x_3|,|x_4|).\]
  $T$ is a $(3,\infty)$ operator because. for any $x$, \[\max\{\norm
  x,\norm{T^2x}\}=\max\{\norm{Tx},\norm{T^3x}\}=
  2\max\{|x_1|,|x_2|,|x_3|,|x_4|\}.\]  
  However, $S=T^2$ is not a $(3,\infty)$ operator because, for
  $x=(0,1,0,0)$ say, we have $\max\{ \norm x, \norm{S^2x}\}=\norm x =
  2$ while $\max\{\norm{Sx},\norm{S^3x}\}=\norm{Sx}=1$.

  Notice that the specific $T$ and $x$ here realises the binary
  $(3,\infty)$ sequence $(1,1,0,0,1,1,0,0,\ldots)
  =(1100)$ while $S=T^2$ has the regular subsequence
  $(10)$. 
  
\end{Proof}

Remark that above $T^3$ is a $(3,\infty)$-isometry.  Indeed, the inverse of an
invertible $(m,\infty)$-isometry must also be $(m,\infty)$.

The definition of a $(3,\infty)$ (binary) sequence may be equivalently
expressed as the interleaving of two binary sequences, neither of
which has consecutive zeros.  

The  graph in Figure~\ref{fig:d1} has a Hamiltonian cycle $P_{(1234567981)}$ which
gives the periodic $(3,\infty)$ sequence $(110011011110)$.

\subsection{$(4,\infty)$ binary sequences}

Any consecutive five terms of a $(4,\infty)$ sequence must have one of
the patterns:
\begin{multicols}{3}
  \begin{enumerate}[$P_\bgroup1\egroup$:\ ]\addtocounter{enumi}{-1}
  \item 0 0 0 1 1
  \item 0 0 1 1 0
  \item 0 0 1 1 1
  \item 0 1 0 0 1
  \item 0 1 0 1 1
  \item 0 1 1 0 0
  \item 0 1 1 0 1
  \item 0 1 1 1 0
  \item 0 1 1 1 1
  \item 1 0 0 1 0
  \item 1 0 0 1 1
  \item 1 0 1 1 0
  \item 1 0 1 1 1
  \item 1 1 0 0 0
  \item 1 1 0 0 1
  \item 1 1 0 1 0
  \item 1 1 0 1 1
  \item 1 1 1 0 0
  \item 1 1 1 0 1
  \item 1 1 1 1 0
  \item 1 1 1 1 1
  \end{enumerate}
\end{multicols}

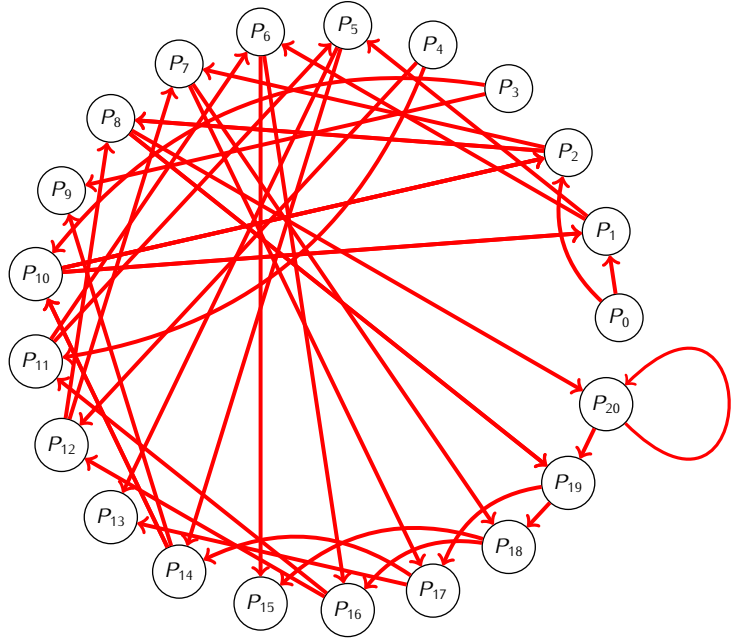
\begin{figure}[h]
 \centering
\begin{tikzpicture}[scale=0.65,transform shape]
   \SetUpEdge[lw         = 1.5pt,
            color      = red,
            labelcolor = white, labelstyle={above}]
            \SetVertexMath
\GraphInit[vstyle=Normal] 
\SetGraphUnit{2}
\tikzset{VertexStyle/.append  style={fill}}
\grEmptyCycle[Math,prefix=P,RA=6,RS=1]{21}
\tikzset{EdgeStyle/.style={->}}
\Edges(P0,P1,P5,P13)
\Edges[style={bend left}](P0,P2)
\Edges(P2,P7,P17,P13)
\Edges(P1,P6,P15)
\Edges(P5,P14,P9)
\Edges(P14,P10,P1)
\Edges(P10,P2,P8,P19)
\Edges[style={bend right}](P19,P17,P14)
\Edges(P14,P10,P1)
\Edges(P8,P20,P19)
\Edges(P19,P18)
\Edges[style={bend right}](P18,P15)
\Edges(P10,P2,P8,P19)
\Edges[style={bend right}](P18,P16)
\Edges(P16,P11,P5)
\Edges(P16,P12,P7,P18)
\Edges(P11,P6,P16)
\Edges(P12,P8)
            \Edges(P3,P9)
            \Edges[style={bend right}](P3,P10)
            \Edges(P10,P2)
            \Edges[style={bend left}](P4,P11)
            \Edges(P4,P12)
   \Loop[dir=EA,style={red,very thick}](P20)
 \end{tikzpicture}
\caption{Graph for $(4,\infty)$ sequence}\label{fig:d2}
\end{figure}

A $(4,\infty)$ sequence is then represented by an infinite path
through the graph shown in Figure~\ref{fig:d2}.

As an example of a cycle, we have the path $P_1\to
P_5\to P_{14}\to P_{10}\to P_1$ which gives the sequence
$(001100110011\ldots)$ or $(0011)$.   Again, the 2-regular subsequence
$(01)$ is not a $(3,\infty)$ sequence.  The cycle $P_2\to P_7\to P_{17}\to
P_{14}\to P_{10} \to P_2$ gives the sequence $(00111)$ 

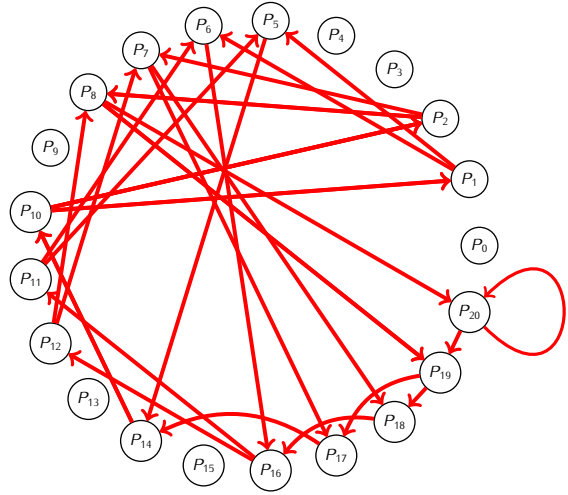
\begin{figure}[t]
 \centering
\begin{tikzpicture}[scale=0.5,transform shape]
   \SetUpEdge[lw         = 1.5pt,
            color      = red,
            labelcolor = white, labelstyle={above}]
            \SetVertexMath
\GraphInit[vstyle=Normal] 
\SetGraphUnit{2}
\tikzset{VertexStyle/.append  style={fill}}
\grEmptyCycle[Math,prefix=P,RA=6,RS=1]{21}
\tikzset{EdgeStyle/.style={->}}
\Edges(P1,P5)
\Edges(P2,P7,P17)
\Edges(P1,P6)
\Edges(P5,P14)
\Edges(P14,P10,P1)
\Edges(P10,P2,P8,P19)
\Edges[style={bend right}](P19,P17,P14)
\Edges(P14,P10,P1)
\Edges(P8,P20,P19)
\Edges(P19,P18)
\Edges(P10,P2,P8,P19)
\Edges[style={bend right}](P18,P16)
\Edges(P16,P11,P5)
\Edges(P16,P12,P7,P18)
\Edges(P11,P6,P16)
\Edges(P12,P8)
            \Edges(P10,P2)
   \Loop[dir=EA,style={red,very thick}](P20)
 \end{tikzpicture}
\caption{Graph for $(4,\infty)$ sequence, degree 2 vertices removed}\label{fig:d3}
\end{figure}

In general, vertices have degree 4, with in-degree and out-degree
both equal to 2.  For some vertices, the degree is 2: $P_0$, $P_3$ and
$P_4$ have in-degree zero so can only appear at the start of the
sequence, while $P_9,P_{13}$ and $P_{15}$ have out-degree zero and cannot appear in a
non-terminating sequence.   If we remove these six ``boundary''
vertices then we obtain the graph in Figure~\ref{fig:d3}.  It is not
Hamiltonian ($P_{10}$ is the unique predecessor of both $P_1$ and $P_2$).   The possible vertex patterns all come about as
the interleaving of two binary sequences which have no consecutive
zeros.  In other words, a (non-terminating) $(4,\infty)$ sequence is a $(3,\infty)$
sequence if it can be extended to the left (i.e. is the $\N$-indexed
subsequence of a $\Z$-indexed $(4,\infty)$-sequence).   This is not
surprising given \cite[Proposition 6.4]{MR2852193} which states that
an invertible $(m,\infty)$-isometry is an $(m-1,\infty)$-isometry if
$m$ is even.

The number of $n$-bit binary numbers without consecutive zeros is
$F_n$, the $n$th Fibonacci number (starting with $F_{-1}=F_{0}=1$).
Consequently, the number of  $2n$-bit (resp. $2n-1)$) numbers which form a (finite)
$(3,\infty)$ sequence is $F_n^2$ (resp. $F_nF_{n-1}$).   Accordingly, the nine
vertices in Figure~\ref{fig:d1}, each representing a 4-bit
$(3,\infty)$-sequence, reflect the fact that $F_2^2=9$. 
Each edge in the graph represents a five bit $(3,\infty)$ sequence and
so  $F_3F_2=15$ is reflected in the 15 edges.  Paths of length $k$
give distinct $(4+k)$-bit $(3,\infty)$ sequences allowing us to say
that the graph contains $F_{2+\lfloor \frac {k+1}2 \rfloor} F_{2+
  \lfloor \frac k2 \rfloor}$ distinct paths of length $k$.

Similar statements can be made for finite $(2m,\infty)$-sequences
(which contain $(2m-1,\infty)$-sequences as a subset).  These are
formed as the interleaving of binary sequences which do not contain
$m$ consecutive zeros.  If $G_n$ denotes the number of $n$-bit binary
sequences which do not contain 3 consecutive zeros then $G_n$
satisfies the generalised Fibonacci relation
$G_{n+1}=G_n+G_{n-1}+G_{n-2}$ starting with $G_0=G_{-1}=1, G_{-2}=0$.
The number of $2n$-bit $(3,\infty)$-sequences is $G_n^2$.
Construction of $(m,\infty)$-operators which reproduce these finite
$(m,\infty)$-sequences or their infinite periodic extensions can be done as in Corollary~\ref{cor:fail}.

Finally, we remark that any subsequence, and in particular any regular subsequence, of the $(m,\infty)$-sequence
\[ (\underlabel{\mbox{\scriptsize $m-1$
      zeros}}{0,0,\ldots,0},1,1,1,\ldots) \]
is clearly also an $(m,\infty)$-sequence.  This expressly provides a
proof of \cite[Proposition 3.10]{MR3947303}:

\begin{Proposition}
  Let $T\in L(X)$ and $m\in\N$, with $m\ge2$ such that \[ \norm{T^m
      x}=\norm{T^{m-1}x} \mbox{ and } \norm{T^mx}\ge \norm{T^\ell x}, \]
    for any $\ell\in \{0,1,\ldots,m-2\}$ and $x\in X$.  Then $T^k$ is an
  $(m,\infty)$-isometry for any $k\in\N$.
\end{Proposition}

Furthermore, with this perspective we may strengthen the conclusion
to $T^k$ is an $(\lceil \frac mk\rceil,\infty)$-isometry.

\def\cprime{$'$}

\bigskip

\end{document}